\documentclass[preprint,12pt,number]{elsarticle}

\usepackage{amsmath} 
\usepackage{amssymb}  
\usepackage{amsfonts}
\usepackage{amsthm}
\usepackage{pifont}
\newtheoremstyle{indented}
    {3pt}
    {3pt}
    {\addtolength{\leftskip}{0em}}
    {}
    {\bfseries}
    {.}
    {1em}
    {}
\theoremstyle{indented}

\newtheorem{theorem}{Theorem}[section]
\newtheorem{corollary}{Corollary}[section]
\newtheorem{proposition}{Proposition}[section]

\newtheorem{definition}{Definition}[section]
\newtheorem{exmp}{Example}[section]

\newtheorem{lemma}{Lemma}[section]

\usepackage{float} 
\usepackage[skip=1pt,font=footnotesize]{caption}

\usepackage{subcaption}

\usepackage{etaremune}

\usepackage{tikz}
\usetikzlibrary{shapes,arrows,calc,positioning}
\tikzstyle{bigblock} = [draw, fill=blue!20, rectangle, 
    minimum height=6em, minimum width=8em]
\tikzstyle{medblock} = [draw, fill=blue!20, rectangle, 
    minimum height=4em, minimum width=4em]    
\tikzstyle{mux} = [draw, fill=black!20, rectangle, 
    minimum height=5em, minimum width=0.1em]    
\tikzstyle{smallblock} = [draw, fill=blue!20, rectangle, 
    minimum height=2em, minimum width=3em]
    
\tikzstyle{data_block} = [draw, fill=green!20, rectangle, 
    minimum height=2em, minimum width=3em]
\tikzstyle{ops_block} = [draw, fill=blue!20, rectangle, 
    minimum height=2em, minimum width=3em]    
\tikzstyle{est_block} = [draw, fill=red!20, rectangle, 
    minimum height=2em, minimum width=3em]    
    
\tikzstyle{sum} = [draw, fill=blue!20, circle, node distance=1cm,minimum height=0.5cm]
\tikzstyle{signal} = [coordinate]
\tikzstyle{pinstyle} = [pin edge={to-,thin,black}]
\tikzstyle{block} = [draw, fill=blue!20, rectangle, 
    minimum height=3em, minimum width=9em]
\tikzstyle{blockS} = [draw, fill=blue!20, rectangle, 
    minimum height=3em, minimum width=4em]    
\tikzstyle{input} = [coordinate]
\tikzstyle{output} = [coordinate]
\usetikzlibrary{matrix}

\usetikzlibrary{positioning}
\usetikzlibrary{math}

\usepackage{hyperref}
\usepackage{xcolor}
\hypersetup{
    colorlinks,
    linkcolor={blue!100!black},
    citecolor={blue!50!black},
    urlcolor={blue!80!black}
}

\newcommand{\bc}{\begin{center}}
\newcommand{\ec}{\end{center}}
\newcommand{\benum}{\begin{enumerate}}
\newcommand{\eenum}{\end{enumerate}}
\newcommand{\nn}{\nonumber}
\newcommand{\matl}{\left[ \begin{array}}
\newcommand{\matr}{\end{array} \right]}
\newcommand{\matls}{\left[ \begin{smallmatrix}}
\newcommand{\matrs}{\end{smallmatrix} \right]}
\newcommand{\isdef}{\stackrel{\triangle}{=}}

\newcommand{\rmT}{{\rm T}}

\newcommand{\BBC}{{\mathbb C}}

\newcommand{\BBN}{{\mathbb N}}

\newcommand{\BBR}{{\mathbb R}}

\newcommand{\SA}{{\mathcal A}}
\newcommand{\SB}{{\mathcal B}}
\newcommand{\SC}{{\mathcal C}}

\newcommand{\SN}{{\mathcal N}}

\newcommand{\SR}{{\mathcal R}}

\renewcommand{\matl}{\begin{bmatrix}}
\renewcommand{\matr}{\end{bmatrix} }

\newcommand{\rank}{{\rm rank  } }

\newcommand{\realization}[4]{
\left[\begin{array}{c|c}
    #1 & #2\\
    \hline
    #3 & #4 \end{array}\right]
}

\newcommand{\EndProofInEq}{\tag*{\mbox{\qed}}}

\newcommand{\wrapmat}[1]{\matl #1 \matr}

\usepackage{comment}
\usepackage{marginnote}

\usepackage{enumitem,amssymb}
\newlist{todolist}{itemize}{2}
\setlist[todolist]{label=$\square$}
\usepackage{pifont}
\usepackage[colorinlistoftodos,textwidth=1 in,textsize=footnotesize]{todonotes}

\usepackage{graphicx}      
\usepackage{textcomp}
\usepackage{calc}
\usepackage{mathtools}
\usepackage{xparse}
\usepackage{amssymb}
\usepackage{stmaryrd}
\usepackage{amsmath}
\usepackage{cases}
\usepackage{mathtools}
\usepackage{cuted}
\usepackage{flushend}
\usepackage{makecell}
\usepackage{booktabs}
\usepackage{array}
\usepackage{arydshln}
\usepackage{siunitx}
\usepackage{xcolor}
\usepackage{caption}

\usepackage[linesnumbered,ruled,vlined]{algorithm2e}
\usepackage{wrapfig}

\SetKwInput{KwInput}{Input}                
\SetKwInput{KwInitialize}{Initialize}                
\SetKwInput{KwOutput}{Output}              

\newcommand{\bluetext}[1]{\textcolor{black}{#1}}

\begin{document}

\begin{frontmatter}


\title{\LARGE \bf
Computing Invariant Zeros of a MIMO Linear System \\
Using State-Space Realization
}


\author{\large Jhon Manuel Portella Delgado and Ankit Goel
\thanks{Jhon Manuel Portella and Ankit Goel are with the Department of Mechanical Engineering, University of Maryland, Baltimore County, MD 21250.
{\tt \small \{jportella,ankgoel\}@umbc.edu}}}


\begin{abstract}
%
%
%
Poles of a multi-input, multi-output (MIMO) linear system can be computed by solving an eigenvalue problem, however, the problem of computing its invariant zeros is equivalent to a generalized eigenvalue problem. 
This paper revisits the problem of computing the invariant zeros by solving an eigenvalue problem. 
We introduce a realization called the \textit{invariant zero form} in which the system's invariant zeros are isolated in a partition of the transformed dynamics matrix.
It is shown that the invariant zeros are then the eigenvalues of a partition of the transformed dynamics matrix.
Although the paper's main result is proved only for square MIMO systems, the technique can be heuristically extended to nonsquare MIMO systems, as shown in the numerical examples. 
%
%
%
%
%
%
\end{abstract}


\end{frontmatter}



















\section{INTRODUCTION}

Zeros are fundamental in the study of systems and control theory. 
While the poles affect the system stability, transients, and convergence rate, zeros affect the undershoot, overshoot, and zero crossings \cite{macfarlane1976poles,desoer1974zeros,tokarzewski2006finite}.
Furthermore, non-minimum phase zeros, which are zeros in the open-right-half-plane, limit performance and bandwidth due to limited gain margin, and exacerbate the tradeoff between the robustness and achievable performance of a feedback control system \cite{hoagg2007nonminimum,havre2001achievable,wonham1970decoupling}.  
Precise knowledge of zeros is thus crucial in the design of reliable control and estimation systems. 

In a single-input, single-output (SISO) linear system, zeros and poles can be computed by computing the roots of the numerator and the denominator of its transfer function, respectively. 
In multi-input, multi-output (MIMO) linear systems, zero computation is a considerably more difficult problem. 
In fact, in MIMO systems, various types of zeros can be defined, such as, 
invariant zeros\cite{owens1977invariant}, 
transmission zeros\cite{davison1974properties}, 
decoupling zeros \cite{horan1980decoupling}, and 
blocking zeros\cite{patel1986blocking}. 
This paper focuses on the computation of invariant zeros from the state-space realization of a linear system. 
%

In a state-space realization of a MIMO system, the poles are the subset of the eigenvalues of the dynamics matrix. Therefore, the computation of poles is an eigenvalue problem.
On the other hand, invariant zeros of a MIMO system with a state-space realization $(A,B,C,D)$ are the complex numbers for which the rank of the Rosenbrock system matrix  
\begin{align}
    \chi(\lambda) 
        \isdef
            \matl 
                \lambda I - A & B \\C & -D
            \matr
\end{align}
drops \cite{laub1978calculation,tokarzewski2011invariant,tokarzewski2009zeros}.  
%
Note that $\chi(\lambda)$ is a pencil
\begin{align}
    \chi(\lambda)
        =
            \lambda
            \matl 
                I & 0\\
                0 & 0
            \matr
            -
            \matl 
                A & -B\\
                -C & D
            \matr,
\end{align}
and thus the invariant zeros of $(A,B,C,D)$ are the 
generalized eigenvalues of the pencil $\chi(\lambda)$
\cite[pp. 1278]{bernstein2009matrix}.
The generalized eigenvalues of the pencil $\chi(\lambda)$ can be computed by decomposing it into the generalized Schur form as shown in \cite{emami1982computation, misra1994computation,misra1989computation,misra1989minimal}. 
However, this approach yields extraneous zeros, which are removed
heuristically \bluetext{by explicitly computing that the rank of $\chi(\lambda)$ at all the resulting candidate zeros and retaining only those at which the rank of $\chi(\lambda)$ drops.}     
\bluetext{An alternative approach, based on the elementary root locus analysis, uses the fact that the closed-loop poles approach the zeros of the system under high-gain output feedback.}
Thus, the zeros of $(A,B,C,D)$ are a subset of the eigenvalues of
$\lim_{K \to \infty } A+ BK(I- D K)^{-1} C$ \cite{davison1978algorithm,garbow1977matrix}.

A geometric approach to compute the invariant zeros of the system by solving an eigenvalue problem is described in \cite{aling1984nine,basile1994controlled,morris2010invariant}.
This paper revisits the problem of invariant zeros computation as an eigenvalue problem.
In contrast to the geometric approach of \cite{aling1984nine}, the {significantly more straightforward algebraic} approach presented in this paper is motivated by and closely related to the normal form of a dynamic system \cite{isidori1985nonlinear}.
The algorithm based on the geometric approach described in \cite{basile1994controlled} and the algorithm presented in this paper to compute the invariant zeros are summarized in Appendix \ref{appndx:Algorithms}.

\bluetext{The normal form of a dynamic system partitions the state into the well-known zero dynamics and a completely linearizable dynamics.
Linearizable dynamics can then be used to construct output tracking controllers as shown in \cite{zhao2024weak, zhou2024invertibility, portella2024circumventing, zhou2023relative, korovin2007canonical}.
%
In our preliminary work, described in \cite{portella2024computing}, it is shown that the partition corresponding to the zero dynamics of a SISO system, in fact, contains the zeros of the system. 
This paper extends this result to MIMO systems.
In particular, we show that the partition corresponding to the MIMO system's zero dynamics contains the system's invariant zeros.
Although this fact is well known, \cite{zhou2023relative,Khalil:1173048}, this work presents an algebra-based proof of this result with minimum assumptions on the structure of the system matrices and extends the application of the result to the nonsquare MIMO systems as well. 
}


%

%
%
%


To compute the invariant zeros, we construct a change of basis that isolates the zeros in a partition of the dynamics matrix in the new basis and show that the invariant zeros of the system are precisely the eigenvalues of this partition. 
Since the invariant zeros of the system are isolated in a partition of the dynamics matrix in the new basis, we call this realization the \textit{invariant zero form}.
The main result of this paper is proved only for square MIMO  systems, however, as shown in Section \ref{sec:examples}, the result can be heuristically extended to nonsquare MIMO systems. 
Finally, since the paper's main result does not require the minimality of the realization, the computed zeros are the invariant zeros of the state-space realization.

%
%
%


%

%
%

The paper is organized as follows. 
Section \ref{sec:notation} presents the notation used in this paper, 
Section \ref{sec:ZSF} introduces the invariant zero decomposition of a MIMO linear system and describes several properties of the resulting realization, 
Section \ref{sec:main_result} presents and proves the main result of the paper, 
and 
Section \ref{sec:examples} presents numerical examples to confirm the main result of this paper.
Finally, the paper concludes with a discussion in Section \ref{sec:conclusion}.

\section{Notation}
\label{sec:notation}
Let $A \in \BBR^{n\times m}.$
Then, $A_{[i,j]}$ is the matrix obtained by removing the $i$th row and $j$th column of $A.$
Note that if $i=0,$ then only the $j$th column is removed, if $j=0,$ then only the $i$th row is removed.
Therefore, $A_{[0,0]} = A.$
$A_{\{i:j\}}$ is the matrix composed of columns $\{i,i+1, \ldots, j\} $ of $A.$
\bluetext{
The set of integers between $n$ and $m$, where $n\leq m,$ that is, $\{n, n+1, \ldots, m\},$ is denoted by
    $\{ n, \ldots, m\}.$}
%
$0_{n \times m}$ denotes the $n \times m$ zero matrix and 
$I_n$ denotes the $n\times n$ identity matrix. 
\bluetext{
$\SN(A)$ and $\SR(A)$ denotes the nullspace and the range space of $A,$ respectively. 
%
${\rm blkdiag}(A_1,A_2,\cdots,A_n)$ denotes a block diagonal matrix formed by aligning $A_1,A_2,\cdots,A_n$ along the diagonal.
${\rm mspec}(A)$ denotes the multiset of eigenvalues of $A.$ 
$A^+$ denotes the Moore-Penrose pseudoinverse of $A$. 
}

\section{Invariant Zero Decomposition}
\label{sec:ZSF}
In this section, we introduce the \textit{invariant zero decomposition}
that transforms a MIMO linear system to its \textit{invariant zero form}
.
%
%
Like the Kalman decomposition that isolates the controllable and observable portion of the dynamics matrix in a partition of the transformed dynamics matrix, 
the invariant zero decomposition isolates the zeros of the system in a partition of the transformed dynamics matrix. 
The invariant zeros of the system can then be computed by solving an eigenvalue problem of a partition of the dynamics matrix expressed in the invariant zero form. 

To facilitate the following discussion, we first precisely define the invariant zeros of a MIMO system. 
The following definition appears in \cite[p.~830]{bernstein2009matrix}.
\begin{definition}
    \label{def:IZ}
    Let $G \in \BBR(s)^{l_y \times l_u}$, where $G \sim 
        \left[
        \begin{array}{c|c}
            A & B \\
            \hline
            C & D
        \end{array} 
        \right]
        .$
    $z \in \mathbb{C}$ is an \textbf{invariant zero} of  $G$ if 
    \begin{align}
        \rank \  \chi(z) < \rank \  \chi,
    \end{align}
    where $\chi(s)$ is the Rosenbrock system matrix 
    \begin{align}
        \chi(s)
            \isdef 
                \matl
                    sI-A & B\\
                    C & -D
                \matr,
    \end{align}
    and $\rank \ \chi \isdef  \underset{s\in \BBC} {\rm max}  \ \rank( \chi(s)).$
    The set of invariant zeros of $G$ is denoted by ${\rm izeros} 
                    \left[
                    \begin{array}{c|c}
                        A & B \\
                        \hline
                        C & D
                    \end{array} 
                    \right].$

\end{definition}

In the following, we construct a state-transformation matrix that transforms a realization of a linear system to its invariant zero form.
%
%
Consider the linear system
\begin{align}
    \lambda x &= A x + Bu,
    \label{eq:state_x} \\
    y &= Cx, \label{eq:output_y}
\end{align}
where $x\in \BBR^{l_x}$ is the state, 
$u \in \BBR^{l_u}$ is the input, 
$y \in \BBR^{l_y}$ is the output, and the operator $\lambda$ is the time-derivative operator or the forward shift operator.
\bluetext{Without loss of generality, in this paper, we assume that $\lambda$ is the time-derivative operator.}
The relative degree of $i$th output $y_i$ is denoted by $\rho_i,$ and the relative degree of $y,$ defined as $\sum_i^{l_y} \rho_i,$ is denoted by $\rho.$
\bluetext{
Note that the relative degree of $y_i$ is defined as the minimum number of times $y_i$ should be differentiated such that one of the inputs appears explicitly in the differential equation \cite{isidori1988nonlinear}. }

%
%
Note that \eqref{eq:state_x}, \eqref{eq:output_y} is a strictly proper system since the feedforward matrix $D = 0.$ 
Although the main result presented in this paper requires $\rho \geq l_y,$ that is, feedforward matrix $D=0,$ this is not a restrictive assumption since, as shown in Appendix \ref{appndx:feedforward}, a dynamic extension of the system by augmenting it with an arbitrary pole renders the feedforward matrix of the augmented system zero without affecting the invariant zeros of the system. 
Note that an invariant zero can therefore be identified despite possible cancellation by the pole used to extend the system.
Example \ref{exmp:D_nonzero} considers the case of an exactly proper system.

Let 
$\overline B \in \BBR^{(l_x -l_u) \times l_x}$ 
be a full rank matrix such that $\overline B B = 0$ 
and let
\begin{align}
    {\overline{C}}
        \isdef 
            \matl
                \overline{C}_1\\
                \overline{C}_2\\
                \vdots\\
                \overline{C}_{l_y}
            \matr \in \BBR^{\rho \times l_x},
        \label{eq:Cbar_def}
\end{align}
where, for 
$ i \in \{ 1, \ldots, l_y\}$
\begin{align}
    \overline{C}_i 
        \isdef 
            \matl
                C_i\\
                C_iA\\
                \vdots\\
                C_iA^{\rho_i-1}
            \matr \in \BBR^{\rho_i \times l_x}.
         \label{eq:Cbar_def_2}
\end{align}
Note that the columns of $B$ belong to $\SN(\overline{B}),$ thus $\overline B$ can be computed using \texttt{null(B')} in MATLAB.
Next, define
\begin{align}
    \overline T 
        \isdef
            \matl
                \overline{B} \\
                \overline{C}
            \matr
            \in \BBR^{(l_x-l_u+\rho) \times l_x}.
    \label{eq:Tbar_def}
\end{align}

\begin{proposition}
    \label{prop:Tbar_full_rank}
    Consider the system \eqref{eq:state_x}, \eqref{eq:output_y}.
    Let $\rho \in [l_y, \infty).$    
    Let $\overline T$ be given by \eqref{eq:Tbar_def}.
    %
%
    Then, 
    \begin{align}
        \rank (\overline T )
            &= 
                l_x - 
                \rank (B) + 
                \rank (\hat C)
                \nn \\ &\quad
                -
                {\rm dim}
                \left( 
                    \SR(\overline{B}^\rmT) \cap 
                    \SR
                    \left(
                \hat C^\rmT
                    \right)
                \right),
        \label{eq:rankTbar}
    \end{align}
    where 
    \begin{align}
        \hat C 
            \isdef 
                \matl 
                    C_1 A^{\rho_1-1} \\
                    \vdots \\
                    C_{l_y} A^{\rho_{l_y}-1} 
                \matr 
            \in \BBR^{l_y \times l_x}.
    \end{align}
\end{proposition}
\begin{proof}
    Note that $\overline B$ has $l_x-\rank(B)$ linearly independent rows. 
    Since the relative degree of $y_i$ is $\rho_i>0,$ it follows that, for each 
        $ j \in \{ 1, \ldots, l_y\},
        i \in \{ 1, \ldots, \rho_j-1\}, $
    $C_j A^{i-1} B = 0,$ which implies that each
    $C_j A^{i-1}$
    is in the row range space of $\overline B$ since columns of $B$ are in $\SN(\overline{B}).$
    Therefore,
    \renewcommand{\arraystretch}{1.5}
    \begin{align}
        \rank 
            \matl
                \overline{B} \\
                \overline{C}
            \matr
                =
                    \rank 
                    \matl
                        \overline{B} \\
                        \hat C
                    \matr,  \nn
    \end{align}
    \renewcommand{\arraystretch}{1}
    and thus \eqref{eq:rankAC} in Lemma \ref{lem:rank_tian} yields \eqref{eq:rankTbar}.
\end{proof}

Next, to construct a state-transformation matrix from the rows of $\overline{T}$, we will restrict the discussion to the cases where $\rank(\overline{T}) = l_x.$
To do so, assume $\overline C$ is full-rank, $\rho \leq l_x, $ and define 
\begin{align}
    T 
        \isdef 
            \matl 
                B_z \\
                \overline C
            \matr 
            \in \BBR^{l_x \times l_x},
    \label{eq:T_def}
\end{align}
where 
$B_z \in \BBR^{l_{\eta} \times l_x},$
where $l_{\eta} \isdef l_x - \rho,$ 
is chosen such that rows of $B_z$ and rows of $\overline C$ are linearly independent. Consequently, $T$ is full-rank and thus invertible.
Note that $\SR(B_z^\rmT) \subseteq \SR(\overline B^\rmT).$

Next, define $\eta \in \mathbb{R}^{l_{\eta}}$ as the first $l_{\eta}$ components of $Tx$ and $\xi \in \mathbb{R}^\rho$ as the rest of $Tx$, that is,
    \begin{align}
        \matl 
            \eta \\
            \xi 
        \matr 
            =
                T x.
        \label{eq:state_transformation}
    \end{align}
    Substituting \eqref{eq:state_transformation} in \eqref{eq:state_x}, \eqref{eq:output_y} yields
    \begin{align}
        \matl 
            \dot \eta \\
            \dot \xi 
        \matr 
            &= 
                \SA 
                \matl 
                    \eta \\
                    \xi 
                \matr 
                + \SB u, \label{eq:state_chi} \\
        y 
            &=
                \SC \matl 
                        \eta \\
                        \xi 
                    \matr  , \label{eq:output_y_chi}
    \end{align}
where
\begin{align}
    \SA \isdef T A T^{-1}, \quad 
    \SB 
        \isdef
            TB, \quad 
    \SC 
        \isdef 
            C T^{-1}.
    \label{eq:ZSF_SS}
\end{align}
Using the transformation matrix $T$, the realization $(A,B,C)$ is transformed to its invariant zero form $(\SA,\SB,\SC)$. 
\bluetext{Note that the variables $\eta$ and $\xi$ have been classically used in the literature to denote the zero dynamics and the linearizable dynamics \cite{isidori1985nonlinear, zhou2023relative}. }

As will be shown in the next section, the invariant zeros of the system are the eigenvalues of the partition of $\SA$ corresponding to $\eta$ dynamics.
Since $\eta = B_z x$ and thus $\eta \in \SR(B_z),$ we call $\SR(B_z)$ the \textit{invariant zero subspace} of the system. 

Next, partition $\SA$ as
\begin{align}
    \SA 
        =
            \matl 
                \SA_{\eta} & \SA_{\eta \xi} \\
                \SA_{\xi \eta} & \SA_{\xi}
            \matr,
\end{align}
where $\SA_\eta$ is the $l_{\eta} \times l_{\eta}$ upper-left block,
$\SA_{\eta \xi}$ is the $l_{\eta} \times \rho$ upper-right block,
$\SA_{\xi \eta}$ is the $\rho \times l_{\eta}$ lower-left block, and 
$\SA_\xi$ is the $\rho \times \rho$ lower-right block of $\SA.$
Next, define $S \isdef T^{-1}$ and partition $S$ as $S = \matl S_\eta & S_\xi \matr,$ where $S_\eta \in \BBR^{l_x \times l_{\eta}}$ contains the first $l_{\eta}\color{black}$ columns of $S$ and $S_\xi \in \BBR^{l_x \times \rho}$ contains the last $\rho$ columns of $S.$
Substituting $T$ and $S$ in \eqref{eq:ZSF_SS} thus implies that
\begin{align}
    \SA_\eta  &=  B_z A S_\eta, \\ 
    \SA_{\eta \xi}  &= B_z A S_\xi, \\ 
    \SA_{\xi \eta}  &= \overline C A S_\eta, \label{eq:SA_xi_eta_def} \\ 
    \SA_\xi  &= \overline C A S_\xi. \label{eq:SA_xi_def} 
\end{align}

The next two results show that $\SA_{\xi\eta},$ $\SA_\xi,$ $\SB,$ and $\SC$ are sparse matrices. 

\begin{proposition}
    \label{prop:SA_xi_eta_form}
    Let $\SA_{\xi \eta}$ and $\SA_\xi$ be defined by \eqref{eq:SA_xi_eta_def} and \eqref{eq:SA_xi_def}.
    Then, 
    \begin{align}        
       \SA_{\xi  \eta}
            &=
                \matl
                    \matl
                        0_{(\rho_1-1) \times l_{\eta}}\\
                        C_1A^{\rho_1} S_\eta
                    \matr
                    \\
                    \vdots \\
                    \matl
                        0_{(\rho_{l_y}-1) \times l_{\eta}}\\
                        C_{l_y}A^{\rho_{l_y}} S_\eta
                    \matr
                \matr
                \in \BBR^{\rho \times l_{\eta}}
                ,
        \label{eq:SA_xi_eta_def_2}
    \\ 
       \SA_{\xi  }
            &=
                \matl
                    \matl 
                    \matl 
                         0_{(\rho_1-1) \times 1} & I_{\rho_1-1} & 0_{(\rho_1-1) \times (\rho-\rho_1)}
                    \matr \\
                    C_1A^{\rho_1} S_\xi
                    \matr
                    \\
                    \vdots\\
                    \matl 
                        \matl 
                            0_{(\rho_{l_y}-1) \times (\rho-\rho_{l_y}+1)} & I_{\rho_{l_y}-1}
                        \matr \\
                        C_{l_y}A^{\rho_{l_y}} S_\xi
                    \matr 
                \matr
    \in \BBR^{\rho \times \rho}.
    \label{eq:SA_xi_def_2}
    \end{align}
\end{proposition}
\textbf{Proof.}
Note that 
\begin{align}
    \matl B_z \\ \overline C \matr \matl S_\eta & S_\xi \matr 
        =
            I_{l_x}
        =
            \matl 
                I_{l_{\eta}} & 0_{l_{\eta} \times \rho } \\
                0_{ \rho \times l_{\eta} } & I_{\rho} 
            \matr, 
    \nn
\end{align}
and thus
%
\begin{align}
    \overline C S_\eta         
        &=
            \matl
                \matl 
                    C_1\\
                    \vdots\\
                    C_1A^{\rho_1-1}
                \matr
                \\
            \vdots\\
                \matl
                    C_{l_y}\\
                    \vdots\\
                    C_{l_y}A^{\rho_{l_y}-1}
                \matr 
        \matr S_\eta=
            0_{\rho \times l_{\eta}},\nn\\ 
    \overline C S_\xi         
        &=
            \matl
                \matl
                    C_1\\
                    \vdots\\
                    C_1A^{\rho_1-1}
                \matr
                \\
            \vdots\\
                \matl 
                    C_{l_y}\\
                    \vdots\\
                    C_{l_y}A^{\rho_{l_y}-1}
                \matr 
        \matr S_\xi=
            I_{\rho}, \nn
\end{align}
which implies that, for each 
    $ j \in \{ 1, \ldots, l_y\},
    i \in \{ 0, \ldots, \rho_j-1\}$
$C_j  A^{i} S_\eta = 0$
and 
$C_jA^i S_\xi $ is the $(\rho_1+\cdots + \rho_{j-1}+i+1)$th row of $I_\rho.$ 
Therefore, 
\begin{align*}
    \SA_{\xi \eta}  
        &=
            \overline C A S_\eta
        = 
            \matl
                \wrapmat
                {
                C_1A S_\eta\\
                \vdots\\
                C_1A^{\rho_1}S_\eta
                }
                \\
                \vdots
                \\
                \wrapmat
                {
                C_{l_y}AS_\eta\\
                \vdots\\
                C_{l_y}A^{\rho_{l_y}}S_\eta
                }
        \matr             
        =
            \matl
                \wrapmat
                {0_{(\rho_1-1) \times l_{\eta}}\\
                    C_1A^{\rho_1} S_\eta}
                \\
                \vdots
                \\
                \wrapmat{
                    0_{(\rho_{l_y}-1) \times l_{\eta}}\\
                    C_{l_y}A^{\rho_{l_y}} S_\eta}               
                \matr,
        \\
    \SA_{\xi}  
        &=
            \overline C A S_\xi
        = 
            \matl
                \wrapmat
                {C_1A S_\xi\\
                \vdots\\
                C_1A^{\rho_1}S_\xi}
                \\
            \vdots\\
                \wrapmat
                {C_{l_y}AS_\xi\\
                \vdots\\
                C_{l_y}A^{\rho_{l_y}}S_\xi}
        \matr\nn\\             
        &=
            \matl
                \wrapmat
                {\matl 
                         0_{(\rho_1-1) \times 1} & I_{\rho_1-1} & 0_{(\rho_1-1) \times (\rho-\rho_1)}
                    \matr \\
                    C_1A^{\rho_1} S_\xi}
                \\
                \vdots
                \\
                \wrapmat
                {\matl 
                        0_{(\rho_{l_y}-1) \times (\rho-\rho_{l_y}+1)} & I_{\rho_{l_y}-1}
                \matr \\
                C_{l_y}A^{\rho_{l_y}} S_\xi}
            \matr.
    \EndProofInEq
\end{align*}

\begin{proposition}
\label{be_def}
    Let $\SB$ and $\SC$ be defined by \eqref{eq:ZSF_SS}.
    Then, 
    \begin{align}
    \SB 
        =
            \matl
                0_{l_{\eta} \times l_u}\\
                \SB_\xi 
            \matr,
            \quad
    \SC 
        =
            \matl 0_{l_y \times l_{\eta}} & C_{\xi}\matr,
\end{align}
where 
\begin{align}
    \SB_\xi  
        \isdef 
            \matl
                \overline B_1 \\
                \vdots\\
                \overline B_{l_y} \\
            \matr
            \in \BBR^{\rho \times l_u}
            , 
    \quad 
    \SC_\xi
        \isdef 
            \matl
                E_1
                & 
                \cdots 
                & 
                E_{l_y}
            \matr
            \in \BBR^{l_y \times \rho} ,
    \label{eq:Bxi_Cxi_def}
\end{align}
and, for 
    $i \in \{ 1, \ldots, l_y\},$
$\overline B_i \isdef \matl 0_{(\rho_{i}-1) \times l_u} \\ C_i A^{\rho_i -1} B \matr$
and
$E_i \isdef \matl e_i & 0_{l_y \times (\rho_i-1)} \matr,$ where $e_i$ is the $i$th column of $I_{l_y}.$
\end{proposition}
\textbf{Proof.}
Note that, for each
    $j \in \{ 1, \ldots, l_y\},$
    $i \in \{ 0, \ldots, \rho_j-2\},$
$C_j A^i B = 0,$ which implies
{
\begin{align}
    \SB
        =
            TB
        =
            \matl
                B_z B\\
                \overline C B
            \matr 
        &=
            \matl
                0_{l_{\eta} \times l_u}\\
                \matl
                    C_1\\
                    C_1A\\
                    \vdots\\
                    C_1A^{\rho_1-1}
                \matr B\\
                \vdots\\
                \matl
                    C_{l_y}\\
                    C_{l_y}A\\
                    \vdots\\
                    C_{l_y}A^{\rho_{l_y}-1}
                \matr
                B
            \matr 
        =
            \matl
               0_{l_{\eta} \times l_u}\\ 
                \matl 0_{\rho_{1}-1 \times l_u} \\ C_1 A^{\rho_1 -1} B \matr\\
                \vdots\\
                \matl 0_{\rho_{l_y}-1 \times l_u} \\ C_{l_y} A^{\rho_{l_y} -1} B \matr
            \matr\nn\\
        &=
            \matl
                0_{l_{\eta} \times l_u}\\
                \overline B_1 \\
                \vdots\\
                \overline B_{l_y} \\
            \matr
        =
            \matl
                0_{l_{\eta} \times l_u}\\
                \SB_\xi 
            \matr.
            \nn
\end{align}
}
Next, since $C$ is composed of the rows $C_1, \ldots, C_{l_y},$
and 
\begin{align*}
    C 
        =
            \SC T 
        =
            \SC  
            \matl
                B_z \\
                \overline C
            \matr 
        =
            \SC
            \matl
                B_z \\
                \matl
                    C_1\\
                    C_1A\\
                    \vdots\\
                    C_1A^{\rho_1-1}
                \matr
                \\    
                \vdots
                \\
                    \matl 
                    C_{l_y}\\
                    C_{l_y}A\\
                    \vdots\\
                    C_{l_y}A^{\rho_{l_y}-1}
                    \matr
                \matr, 
\end{align*}
it follows that 
\begin{align}
    \SC 
        &= 
            \matl
                0_{l_y \times l_{\eta}}
                &
                e_1 & 0_{l_y \times (\rho_1-1)}
                & 
                \cdots 
                & 
                e_{l_y} & 0_{l_y \times (\rho_{l_y}-1)}
            \matr\nn\\
        &= 
            \matl
            0_{l_y \times l_{\eta}}
            &
            E_1
            & 
            \cdots 
            & 
            E_{l_y}
\matr.
\end{align}
\qed

In summary, the invariant zero decomposition is
\begin{align}
    \dot \eta 
        &=
            A_\eta \eta + A_{\eta \xi} \xi,
    \label{eq:ZSF_eta_dot}
    \\
    \dot \xi 
        &=
           \SA_{\xi  \eta} \eta + A_\xi \xi + \SB_\xi  u, 
    \label{eq:ZSF_xi_dot}
    \\
    y 
        &=
            \SC_\xi \xi.
    \label{eq:ZSF_y}
\end{align}
Furthermore, $\SA_{\xi \eta } , \SA_\xi, \SB_\xi, $ and $\SC_\xi$ have sparse structures given by Propositions  \ref{prop:SA_xi_eta_form} and \ref{be_def}, respectively. 
As shown in the next section, in the case of a square MIMO system, the eigenvalues of $\SA_\eta$ are the invariant zeros of the MIMO system.


\begin{exmp}
    Consider the system
    \begin{align}
        G(s)
            =
                \dfrac{s+1}{(s+2)(s+3)(s+4)}
                .\nn
    \end{align}
    Note that $\rho =2$.
    A realization of $G$ is 
    \begin{align*}
        A
            &=
                \matl
                    -9 & -6.5 & -3\\
                    4 & 0 & 0\\
                    0 & 2 & 0
                \matr, \quad
        B
            =
                \matl
                    0.5\\
                    0\\
                    0
                \matr, 
        \\
        C
            &=
               \matl
                    0 & 0.5 & 0.25
                \matr.
    \end{align*}
    Invariant zero decomposition yields
    \begin{align}
        \SA 
            &=
                \left[\begin{array}{c|cc}
                   -1 & 0 & 2\\
                   \hline
                   0 & 0 & 1\\
                   3 & -24 & -8
                   \end{array}
                \right],
        \quad
        \SB
            =
                \left[\begin{array}{c}
                    0\\
                    \hline\\[-0.3 cm]
                    \begin{array}{c}
                        0\\
                        1
                    \end{array}
                \end{array}
                \right],\nn\\
        \SC
            &=
                \left[\begin{array}{c|cc}
                     0 & 1 & 0
                \end{array}
                \right]
    \end{align}
    Note that $\SA_{\xi \eta } , \SA_\xi, \SB_\xi, $ and $\SC_\xi$ have sparse structures given by Propositions  \ref{prop:SA_xi_eta_form} and \ref{be_def}. 
\end{exmp}

\section{Main Result}\label{sec:main_result}
In this section, we prove the main result of the paper, 
which states that, in the invariant zero form, the invariant zeros of a square MIMO system are the eigenvalues of a partition of its dynamics matrix.

The following three results are used in the proof of Theorem \ref{theorem:zeros_MIMO}.

\begin{proposition}
    \label{prop:B_xi_B_xi_plus_Aen}
    Let $\SB_\xi $, defined by \eqref{eq:Bxi_Cxi_def}, be full column rank and 
    $A_{\xi \eta}$ be given by \eqref{eq:SA_xi_eta_def_2}.
    Let $l_u \geq l_y.$ 
    Then, 
    \begin{align}
        (I-\SB_\xi  \SB_\xi  ^{{+}})
       \SA_{\xi  \eta}
            =
                0.
    \label{eq:B_xi_B_xi_plus_Aen}
    \end{align}
\end{proposition}
\begin{proof}
    Since $\SB_\xi$ is full column rank, 
    \begin{align}
        \SB_\xi  \SB_\xi  ^{{+}}
            =
                \SB_\xi  (\SB_\xi  ^\rmT \SB_\xi  )^{-1}\SB_\xi  ^\rmT 
            =
                \text{blkdiag}(J_1, J_2, \ldots, J_{l_y})
            \in \mathbb{R}^{\rho \times \rho},\nn
    \end{align}
    where, for 
        $i \in \{ 1, \ldots, l_y\}$
    \begin{align}
        J_i 
            \isdef
                \matl 
                    0_{(\rho_i-1) \times (\rho_i-1)} & 0_{(\rho_i-1) \times 1}\\
                    0_{1 \times (\rho_i-1)} & 1\\    
                \matr
                \in \mathbb{R}^{\rho_i \times \rho_i},\nn
    \end{align}
    which implies that 
    $
        I_{\rho}-\SB_\xi  \SB_\xi  ^{{+}} 
            =
                \text{blkdiag}(I_{\rho_1}-J_1, \ldots, I_{\rho_{l_y}} - J_{l_y}).
    $
    %
    Using $\SA_{\xi \eta }$ given by \eqref{eq:SA_xi_eta_def_2} yields \eqref{eq:B_xi_B_xi_plus_Aen}.
\end{proof}

\begin{proposition}
\label{prop:rank_big_term_final}
Let $l_u=l_y$.
For all $s \in \mathbb{C},$. Define
\begin{align}
    \Phi(s)&\isdef\left(I_{l_x}-\matl
             sI_{l_{\eta}}-\SA_{\eta} & 0\\
                    -\SA_{\xi \eta} & \SB_\xi   \matr\matl
             sI_{l_{\eta}}-\SA_{\eta} & 0\\
                    -\SA_{\xi \eta} & \SB_\xi   \matr^{{+}}\right), 
    \label{eq:Phi_def}
    \\
    \Psi(s) &\isdef\matl
            -A_{\eta \xi}\\
            sI-A_{\xi} \matr, 
    \label{eq:Psi_def}
    \\
    \Xi
        &\isdef
            I_{\rho}-\SC_\xi^{{+}}\SC_\xi
            .
    \label{eq:Xi_def}
\end{align}
Then,
{
    \begin{align}
        \rank
        \left[
            \Phi(s)
            \Psi(s)
            \Xi
        \right]
                =
                    \rho-l_y.
    \end{align}
}
\end{proposition}
\textbf{Proof.}
It follows from Proposition \ref{prop:SA_xi_eta_form} that $\Psi(s)$ can be written in the form given by  \eqref{eq:bigXi},
\begin{table*}
\centering
\begin{minipage}{0.75\textwidth}
{
\begin{align}
    \Psi(s) 
        =
            \matl   
            -{\SA_{\eta\xi}}\\
            \matl 
                s\matl
                    I_{\rho_1-1} & 0_{(\rho_1-1)\times (\rho-\rho_1+1)}
                    \matr-\matl 
                         0_{(\rho_1-1) \times 1} & I_{\rho_1-1} & 0_{(\rho_1-1) \times (\rho-\rho_1)}
                    \matr \\
                   s{(e_{\overline{\rho}_1})}^\rmT - C_1A^{\rho_1} S_\xi
            \matr        
            \\
            \vdots
            \\
            \matl 
                s\matl
                    0_{(\rho_{l_y}-1)\times (\rho-\rho_{l_y})} & I_{\rho_{l_y}-1} & 0_{(\rho_{l_y}-1)\times 1}
                \matr - \matl 
                            0_{(\rho_{l_y}-1) \times (\rho-\rho_{l_y}+1)} & I_{\rho_{l_y}-1}
                        \matr \\
                        s{(e_{\overline{\rho}_{l_y}})}^\rmT -C_{l_y}A^{\rho_{l_y}} S_\xi
            \matr
        \matr.
    \label{eq:bigXi}
\end{align}
}
\medskip
\end{minipage}
\end{table*}
where, for each 
    $j \in \{ 1, \ldots, l_y\},$
$\overline{\rho}_{j} \isdef \sum_{j=1}^{l_y} \rho_j.$ 
Note that $\overline{\rho}_{l_y} = \rho.$
%
Next, it follows from Proposition \ref{prop:Aplus_A_zeros_column} that if the $i$th column of $\SC_\xi$ is zero, then the $i$th row and $i$th column of $\SC_\xi ^+ \SC_\xi$ is zero, which implies that the $(i,i)$th entry of $\Xi$ is 1.
%
%

Next, define 
\begin{align}
    M(s) 
        \isdef 
            \matl
                sI_{l_{\eta}}-\SA_{\eta} & 0\\
                -\SA_{\xi \eta} & \SB_\xi   
            \matr
            =    
            \matl   
                \matl 
                     sI_{l_{\eta}}- A_{\eta} &  0_{l_{\eta} \times l_u} 
                \matr 
                \\
                0_{(\rho_1-1) \times (l_{\eta}+l_u)} \\
                \matl 
                     C_1A^{\rho_1} S_\eta &  C_1A^{\rho_1-1} B
                \matr 
                \\
                \vdots
                \\
                0_{(\rho_{l_y}-1) \times (l_{\eta}+l_u)} \\
                \matl 
                     C_{l_y}A^{\rho_{l_y}} S_\eta &  C_{l_y}A^{\rho_{l_y}-1} B
                \matr 
            \matr. \nn
\end{align}
Note that $M(s)$ has at least $\rho-l_y$ zero rows.
Next, note that $\Phi(s) = I_{l_x} - M(s) M(s)^+.$

Consider the case where $s \notin {\rm mspec}(\SA_\eta).$ 
In this case, for all $s\in \SC,$ $sI_{l_{\eta}}- \SA_\eta$ is full rank, and thus none of the first $l_{\eta}$ rows of $M(s)$ are zero. 
Then, it follows from Proposition \ref{prop:A_Aplus_zeros} that 
if the $i$th row $M(s)$ is zero, then the $i$th row and $i$th column of $M(s)M(s)^+$ is zero, and thus the $(i,i)$th entry of $\Phi(s)$ is 1.
Therefore, in the case where  $s \notin {\rm mspec}(\SA_\eta),$ 
    \begin{align}
        \Phi(s)
            =
                \matl
                    0_{l_{\eta} \times l_x}\\
                    \matl
                        N_1 \\
                        0_{1 \times l_x}
                    \matr 
                    \\
                    \vdots\\
                    \matl
                        N_{l_y} \\
                        0_{1 \times l_x}
                    \matr 
                \matr,
        \nn
    \end{align}
where, for each 
    $i \in \{ 1, \ldots, l_y\},$
\begin{align}
    N_i
        \isdef 
            \matl 
                e_{l_{\eta}+\overline{\rho}_{i-1} + 1}^\rmT \\
                \vdots\\
                e_{l_{\eta}+\overline{\rho}_{i}-1}^\rmT 
            \matr. \nn
\end{align}
Using \eqref{eq:bigXi}, it follows that $\Phi(s)\Psi(s)$ is given by \eqref{eq:Phi_Psi}, where $\Delta_n(s)$ is the $n \times n$ matrix given by \eqref{def:delta_n} in Appendix \ref{prop:delta_n_rank}, and thus 
$\Phi(s)\Psi(s)\Xi$ is given by \eqref{non_drop}
%
%
It follows from Proposition \ref{prop:delta_n_rank} that, for each 
    $i \in \{ 1, \ldots,l_y\},$
$\rank(\Delta_{\rho_i-1}(s)) = \rho_i-1,$ and thus the structure of $\Phi(s)\Psi(s) \Xi$ given by \eqref{non_drop} implies that  
\begin{align}
    \rank\left( \Phi(s)\Psi(s)\Xi \right)
        =
            \sum_{i=1}^{l_y} (\rho_i-1)
        =
            \rho - l_y. \nn
\end{align}
\begin{table*}
\centering
\begin{minipage}{0.75\textwidth}
\begin{align}
    \Phi(s)\Psi(s)&=\matl
        0_{l_{\eta} \times \rho}\\
        \matl
        s\matl
                    I_{\rho_1-1} & 0_{(\rho_1-1)\times (\rho-\rho_1+1)}
                    \matr-\matl 
                         0_{(\rho_1-1) \times 1} & I_{\rho_1-1} & 0_{(\rho_1-1) \times (\rho-\rho_1)}
                    \matr\\
        0_{1 \times \rho}
        \matr
        \\
        \vdots
        \\
        \matl 
        s\matl
                0_{(\rho_{l_y}-1)\times (\rho-\rho_{l_y})} & I_{\rho_{l_y}-1} & 0_{(\rho_{l_y}-1)\times 1}
            \matr - \matl 
                        0_{(\rho_{l_y}-1) \times (\rho-\rho_{l_y}+1)} & I_{\rho_{l_y}-1}
                    \matr \\
        0_{1 \times \rho}
        \matr
    \matr.
    \label{eq:Phi_Psi}
\end{align}
\medskip
\end{minipage}
\end{table*}
\begin{table*}
\centering
\begin{minipage}{0.75\textwidth}
\begin{align}
    \Phi(s)\Psi(s) \Xi
        =
            \matl
                0_{lz \times \rho} 
                \\
                \line(1,0){240}\\
                \begin{matrix}
                \matl 
                    0_{\rho \times 1}
                    &\matl
                        \matl 
                            \Delta_{\rho_1-1}(s)\\
                            0_{1 \times (\rho_1-1)}
                        \matr 
                        \\
                        0_{\rho_2 \times (\rho_1-1)}\\
                        \vdots\\
                        0_{\rho_{l_y} \times (\rho_1-1)}
                    \matr
                \matr 
                &\cdots
                &
                \matl 
                    0_{\rho \times 1}
                    &\matl
                        0_{\rho_1 \times (\rho_{l_y}-1)}\\
                        0_{\rho_2 \times (\rho_{l_y}-1)}\\
                        \vdots\\
                        \matl 
                            \Delta_{\rho_{l_y}-1}(s)\\
                            0_{1 \times (\rho_{l_y}-1)}
                        \matr
                    \matr
                \matr
                \end{matrix}  
            \matr \in \mathbb{R}^{l_x \times \rho}.
    \label{non_drop}
\end{align}
\medskip
\end{minipage}
\end{table*}

Next, consider the case where $s\in \text{mspec} (\SA_{\eta})$.
In this case, $sI_{l_{\eta}}- \SA_\eta$ is not full rank.
Consequently, the first $l_{\eta}$ rows of $\Phi(s)$ are not necessarily zero, which implies that $\Phi(s)$ can be written in the form given by
\begin{align}
        \Phi(s)
            =
                \matl
                    \matl 
                    \Upsilon(s) & 0_{l_{\eta} \times \rho}
                \matr\\
                    \matl
                        N_1 \\
                        0_{1 \times l_x}
                    \matr 
                    \\
                    \vdots\\
                    \matl
                        N_{l_y} \\
                        0_{1 \times l_x}
                    \matr 
                \matr,
        \nn
    \end{align}
where $\Upsilon(s)$ is an $l_{\eta} \times l_{\eta}$ matrix.
Using \eqref{eq:bigXi}, it follows that $\Phi(s)\Psi(s)$ is given by \eqref{Phi_Psi_exp2},
\begin{table*}
\centering
\begin{minipage}{0.75\textwidth}
{
\begin{align}
    \Phi(s)\Psi(s)&=\matl
        \Omega(s)\\
        \matl
        s\matl
                    I_{\rho_1-1} & 0_{(\rho_1-1)\times (\rho-\rho_1+1)}
                    \matr-\matl 
                         0_{(\rho_1-1) \times 1} & I_{\rho_1-1} & 0_{(\rho_1-1) \times (\rho-\rho_1)}
                    \matr\\
        0_{1 \times \rho}
        \matr
        \\
        \vdots
        \\
        \matl 
        s\matl
                0_{(\rho_{l_y}-1)\times (\rho-\rho_{l_y})} & I_{\rho_{l_y}-1} & 0_{(\rho_{l_y}-1)\times 1}
            \matr - \matl 
                        0_{(\rho_{l_y}-1) \times (\rho-\rho_{l_y}+1)} & I_{(\rho_{l_y}-1)}
                    \matr \\
        0_{1 \times \rho}
        \matr
    \matr.
    \label{Phi_Psi_exp2}
\end{align}
}
\medskip
\end{minipage}
\end{table*}
where 
$\Omega(s) \isdef  -\Upsilon(s) \SA_{\eta \xi},$
%
%
and 
\begin{table*}
\centering
\begin{minipage}{0.75\textwidth}
{
\begin{align}
    \Phi(s)\Psi(s) \Xi
        =
            \matl
                \matl
                0_{l_x \times 1}
                &\matl
                    \Omega_{ \{2:\rho_1\} }(s)\\
                    \matl 
                        \Delta_{\rho_1-1}(s)\\
                        0_{1 \times (\rho_1-1)}
                    \matr \\
                    0_{\rho_2 \times (\rho_1-1)}\\
                    \vdots\\
                    0_{\rho_{l_y} \times (\rho_1-1)}
                \matr
                \matr
                &\cdots
                \matl
                0_{l_x \times 1}
                &\matl
                    \Omega_{ \{\overline{\rho}_{l_y-1}+2:\overline{\rho}_{l_y} \} } (s)\\
                    0_{\rho_1 \times (\rho_{l_y}-1)}\\
                    0_{\rho_2 \times (\rho_{l_y}-1)}\\
                    \vdots\\
                    \matl 
                        \Delta_{\rho_{l_y}-1}(s)\\
                        0_{1 \times (\rho_{l_y}-1)}
                    \matr
                \matr
                \matr
            \matr \in \mathbb{R}^{l_x \times \rho}.
    \label{drop}
\end{align}
}
\medskip
\end{minipage}
\end{table*}
the structure of $\Phi(s)\Psi(s) \Xi$ given by \eqref{drop} implies that  
\begin{align}
    \rank\left(\Phi(s)\Psi(s) \Xi \right)
        =
            \sum_{i=1}^{l_y} (\rho_i-1)
        =
            \rho - l_y. \nn
\end{align}

\qed

%

\begin{theorem}
\label{theorem:zeros_MIMO}
Consider the square MIMO system \eqref{eq:state_x}, \eqref{eq:output_y}
and its invariant zero decomposition, given by \eqref{eq:ZSF_eta_dot}-\eqref{eq:ZSF_y}, 
Then, the invariant zeros of the system \eqref{eq:state_x}, \eqref{eq:output_y} are the eigenvalues of $\SA_\eta$, that is, 
\begin{align}
        {\rm izeros } \left( \realization{A}{B}{C}{0} \right) = {\rm mspec } (\SA_\eta)\nn. 
    \end{align}
\end{theorem}
\textbf{Proof.}
Since invariant zeros are invariant under state transformation, the realizations $\realization{A}{B}{C}{0}$ and $\realization{\SA}{\SB}{\SC}{0}$ have same invariant zeros.
Using the invariant zero decomposition \eqref{eq:ZSF_eta_dot}-\eqref{eq:ZSF_y}, 
the Rosenbrock system matrix of the realization $\realization{\SA}{\SB}{\SC}{0}$ is
    \begin{align}
        \chi (s) 
            =
                \matl
                    sI_{l_{\eta}}-\SA_{\eta} & -A_{\eta \xi} & 0_{(l_x-\rho) \times l_u}\\
                    -\SA_{\xi \eta} & sI-A_{\xi} & \SB_\xi  \\
                    0_{l_y \times (l_x-\rho)} & C_{\xi} & 0_{l_y \times l_u},
                \matr. \nn
    \end{align}
    First, note that it follows from Lemma \ref{lem:rank_tian} in Appendix \ref{appndx:LinearAlgebraFacts} and Proposition \ref{prop:B_xi_B_xi_plus_Aen} that 
    \begin{align}
        \rank 
        \matl
            -\SA_{\xi \eta} & \SB_\xi  \\
            sI_{l_{\eta}}-\SA_{\eta} & 0
        \matr
            &=
                \rank (sI_{l_{\eta}}-\SA_{\eta}) 
                + \rank (\SB_\xi  )                 
            \nn \\ &\quad 
                + \rank \Big[-[I-\SB_\xi  \SB_\xi  ^{{+}}]A_{\xi \eta}[I
            \nn \\ &\quad 
                - (sI_{l_{\eta}}-\SA_{\eta})^{{+}}(sI_{l_{\eta}}-\SA_{\eta})]\Big]
            \nn \\
            &=
                \rank (sI_{l_{\eta}}-\SA_{\eta}) 
                + \rank (\SB_\xi  ).
        \label{eq:rank_sI-An_term}
    \end{align}
    Next, it follows from \eqref{eq:rank_sI-An_term}, 
    Lemma \ref{lem:rank_tian}, 
    and Proposition \ref{prop:rank_big_term_final} that
    \begin{align}
        \rank~\chi(s)
            &=
                \rank 
                \left[\begin{array}{c|cc}
                     -A_{\eta \xi} & sI_{l_{\eta}}-\SA_{\eta} & 0 \\
                     sI-A_{\xi} & -\SA_{\xi \eta} &  \SB_\xi  \\
                    \hline
                 C_{\xi} & 0 & 0
                \end{array}\right]
            \nn \\
            &=
                \rank
                \matl
                    sI_{l_{\eta}}-\SA_{\eta} & 0\\
                    -\SA_{\xi \eta} & \SB_\xi  
                \matr
                + \rank (\SC_\xi)
            \nn \\ &\quad
                + \rank (\Phi(s) \Psi(s)\Xi)
            \nn \\
            &=
                \rank (sI_{l_{\eta}}-\SA_{\eta}) 
                + \rank (\SB_\xi  ) 
                + l_y 
                + \rho - l_y
            \nn \\
            &=
                \rank (sI_{l_{\eta}}-\SA_{\eta}) 
                + \rank (\SB_\xi  ) 
                + \rho.
        \label{eq:rank_chi_final}
    \end{align}
    Note that the rank of $\SC_\eta$ is $l_y$ by construction. 
    It thus follows from \eqref{eq:rank_chi_final} that $\rank~\chi(s)$ drops if and only if $s$ is an eigenvalue of $\SA_\eta,$ thus implying that the invariant zeros \eqref{eq:state_x}, \eqref{eq:output_y} are the eigenvalues of $\SA_\eta.$ 
    \qed

\begin{corollary}
    \label{cor:Num_of_IZ}
    The number of invariant zeros in a square MIMO system is ${\rm max}(l_x - \rho,0).$
\end{corollary}




The invariant zero decomposition does not exist in the case where $l_u > l_y,$ since the rank of $T,$ given by \eqref{eq:T_def}, is necessarily less than $l_x. $ 
However, by augmenting the rows of $T$ with a sufficient number of linearly independent vectors from the null space of $T,$ a nonsingular transformation matrix may be constructed. 
In this case, the resulting transformation may yield the invariant zeros in the partition of the transformed dynamics matrix. 
On the other hand, in the case where $l_u < l_y,$ rank of $T$ can be $l_x$, and thus the invariant zero form does exist, but the eigenvalues $\SA_\eta$ are not necessarily the invariant zeros.
However, as shown in the next section, an ad-hoc technique that squares a nonsquare MIMO system can be used to compute the invariant zeros of the nonsquare system.

\section{Numerical Examples}
\label{sec:examples}
This section presents several numerical examples to verify the paper's main result and extend it to exactly proper and nonsquare systems. 

\begin{exmp}\textbf{[Square MIMO System.]}
    Consider the square MIMO system
    \begin{align}
        G(s)=\dfrac{64}{s^3 + 24 s^2 + 176 s + 384}\matl
            1 & s+4\\
            s-2 & s-8
        \matr,\nn
    \end{align}
    with two inputs and two outputs.
    Note that $\rho_1 =2$ and $\rho_2 =2,$ and thus $\rho = 4.$
    A realization of $G$, computed with MATLAB's \texttt{ss} routine, is 
    \begin{align*}
        A
            &=
                \matl
                    -24 & -11 & -6 & 0 & 0 & 0\\
                    16 & 0 & 0 & 0 & 0 & 0\\
                    0 & 4 & 0 & 0 & 0 & 0\\
                    0 & 0 & 0 & -24 & -11 & -6\\
                    0 & 0 & 0 & 16 & 0 & 0\\
                    0 & 0 & 0 & 0 & 4 & 0
                \matr, \quad
        B
            =
                \matl
                    1 & 0\\
                    0 & 0\\
                    0 & 0\\
                    0 & 2\\
                    0 & 0\\
                    0 & 0
                \matr, 
        \\
        C
            &=
               \matl
                    0 & 0 & 1 & 0 & 2 & 2\\
                    0 & 4 & -2 & 0 & 2 & -4
                \matr.
    \end{align*}
    Letting 
    \begin{align}
        B_z
            =
                \matl 
                    0 & 0 & 1 & 0 & 0 & 0\\
                    0 & 1 & 0 & 0 & 0 & 0
                \matr,
    \end{align}
    it follows from \eqref{eq:T_def} that 
    \begin{align}
        T
            =
                \matl
                    B_z\\
                    C_1\\
                    C_1A\\
                    C_2\\
                    C_2A
                \matr
            =
                \matl
                    0 & 0 & 1 & 0 & 0 & 0\\
                    0 & 1 & 0 & 0 & 0 & 0\\
                    0 & 0 & 1 & 0 & 2 & 2\\
                    0 & 4 & 0 & 32 & 8 & 0\\
                    0 & 4 & -2 & 0 & 2 & -4\\
                    64 & -8 & 0 & 32 & -16 & 0
                \matr,
            \label{eq:T_example_1}
    \end{align}
    and thus
    \begin{align*}
        \SA
            &=
           \left[ \begin{array}{cc|cccc}
                    0 & 4 & 0 & 0 & 0 & 0\\
                    0 & -1 & 2 & -0.25 & 1 & 0.25\\
                    \hline
                    0 & 0 & 0 & 1 & 0 & 0\\
                    96 & 76 & -88 & -21 & 4 & 1\\
                    0 & 0 & 0 & 0 & 0 & 1\\
                    -288 & -536 & -272 & -6 & -88 & -26
            \end{array}
            \right],
        \nn\\
        \SB
            &=
                \left[
                \begin{array}{cc}
                    0 & 0\\
                    0 & 0\\
                    \hline
                    0 & 0\\
                    0 & 64\\
                    0 & 0\\
                    64 & 64
                \end{array}
                \right],\nn 
                \\
        \SC
            &=
            \left[    \begin{array}{cc|cccc}
                    0 & 0 & 1 & 0 & 0 & 0\\
                    0 & 0 & 0 & 0 & 1 & 0
                \end{array}
                \right].\nn
    \end{align*}
    Note that $\SA_\eta = \matl
        0 & 4\\
        0 & -1
    \matr$ and thus ${\rm mspec}(\SA_\eta) = \{-1,0\}.$ 
    MATLAB's \texttt{tzero} routine yields ${\rm izeros} 
        \left[
        \begin{array}{c|c}
            A & B \\
            \hline
            C & D
        \end{array} 
        \right] = \{-1,0\},$ 
    which confirms Theorem \ref{theorem:zeros_MIMO}.
    
    Furthermore, for $s = 8 \notin {\rm mspec}(\SA_\eta),$
    \begin{align}
        \Phi(s) &= \matl
            \begin{array}{cccccc}            0  &  0 &  0   & 0 & 0 & 0\\
            0 & 0 & 0 & 0 & 0 & 0
            \end{array}\\
            \matl
             0 & 0 & 1 & 0 & 0 & 0\\
             0 & 0 & 0 & 0 & 0 & 0
            \matr\\[2.7mm]
            \matl
             0 & 0 & 0 & 0 & 1 & 0\\
             0 & 0 & 0 & 0 & 0 & 0
            \matr
        \matr,
        \nn \\
        \Xi(s) 
            &=  
               \left[ \begin{array}{cccc}
                \begin{array}{@{}*{5}{S[table-format=3.0]}@{}}
                    0 & 0 & 0 & 0\\
                    -2 & 0.25 & -1 & -0.25
            \end{array}\\
                \left[    \begin{array}{@{}*{5}{S[table-format=3.0]}@{}}
                    8 & -1 & 0 & 0\\
                    88 & 29 & -4 & -1
                    \end{array}\right]\\[2.8 mm]
                \left[    \begin{array}{@{}*{5}{S[table-format=3.0]}@{}}
                    0 & 0 & 8 & -1\\
                    272 & 6 & 88 & 34
                    \end{array}\right]
                \end{array}
                \right],\nn\\
        \Psi &=
            \matl
                 0 & 0 & 0 & 0\\
                 0 & 1 & 0 & 0\\
                 0 & 0 & 0 & 0\\
                 0 & 0 & 0 & 1
            \matr,
    \end{align}
    and thus 
    \begin{align}
        \Phi(s) \Xi(s) \Psi &= 
        \left[
        \begin{array}{ccccc}
        \begin{matrix}
            0 & 0 \\
            0 & 0
        \end{matrix} & \begin{matrix}
            0 & 0 \\
            0 & 0
        \end{matrix}\\
            \hline
         \matl
            0 & -1\\
            0 & 0\\
            0 & 0\\
            0 & 0
        \matr& \matl
            0 & 0\\
            0 & 0\\
            0 & -1\\
            0 & 0
        \matr
        \end{array}
        \right],
    \end{align}
    which confirms \eqref{non_drop}, and thus $\rank (\Phi(s) \Psi(s) \Xi) = 2 = \rho-l_y, $ which confirms Proposition \ref{prop:rank_big_term_final}. 
    Next, for $s = -1 \in {\rm mspec}(\SA_\eta),$
    \begin{align}
        \Phi(s) 
            &=
            \matl \begin{array}{cccccc}
                    0 & 0 & 0 & 0 & 0 & 0\\
                    0 & 1 & 0 & 0 & 0 & 0
            \end{array}\\
            \matl
                 0 & 0 & 1 & 0 & 0 & 0\\
                 0 & 0 & 0 & 0 & 0 & 0
            \matr\\[2.8mm]
            \matl
                 0 & 0 & 0 & 0 & 1 & 0\\
                 0 & 0 & 0 & 0 & 0 & 0
            \matr
            \matr,
        \nn \\
        \Xi(s) 
            &=  
                \left[
               \begin{array}{cccc} \begin{array}{@{}*{5}{S[table-format=3.0]}@{}}
                    0 & 0 & 0 & 0\\
                    -2 & 0.25 & -1 & -0.25
            \end{array}\\
                \left[    \begin{array}{@{}*{5}{S[table-format=3.0]}@{}}
                    -1 & -1 & 0 & 0\\
                    88 & 20 & -4 & -1
                    \end{array}\right]\\[2.8 mm]
                \left[    \begin{array}{@{}*{5}{S[table-format=3.0]}@{}}
                    0 & 0 & -1 & -1\\
                    272 & 6 & 88 & 25
                    \end{array}\right]
                \end{array}
                \right],
        \nn\\ 
        \Psi 
            &=
                \matl
                    0 & 0 & 0 & 0\\
                    0 & 1 & 0 & 0\\
                    0 & 0 & 0 & 0\\
                    0 & 0 & 0 & 1
                \matr,
    \end{align}
    and thus 
    \begin{align}
        \Phi(s) \Xi(s) \Psi 
            &=
                \matl
                    \matl
                        0 & 0\\
                        0 & 0.25\\
                        0 & -1\\
                        0 & 0\\
                        0 & 0\\
                        0 & 0
                    \matr&
                    \matl
                        0 & 0\\
                        0 & -0.25\\
                        0 & 0\\
                        0 & 0\\
                        0 & -1\\
                        0 & 0
                    \matr
                \matr,
    \end{align}
    which confirms \eqref{drop}.
    Note that $\rank (\Phi(s) \Psi(s) \Xi) = 2 = \rho - l_y,$ which confirms Proposition \ref{prop:rank_big_term_final}. 

    \hfill{\huge$\diamond$}.
\end{exmp}

\begin{exmp} 
\textbf{[Square MIMO system with Direct Feedthrough.]}
\label{exmp:D_nonzero}
Consider the system
\begin{align}
    G(s)=\dfrac{64}{s^2 + 12 s + 32}\matl
        s + 2 &  s^2 + 2 s + 4\\
        s^2 + 4 s + 8 &  s + 4
    \matr,
\end{align}
Note that $\rho_1=0, \rho_2=0$, and thus $\rho=0$. A realization of $G$, computed with MATLAB's \texttt{ss} routine, is
\begin{align}
    A&=\matl
       -12 &  -4 & 0 & 0\\
       8 & 0 & 0 & 0\\
       0 & 0 & -12 & -4\\
       0 & 0 & 8 & 0
    \matr, \quad 
    B=\matl
        16 & 0\\
        0 & 0\\
        0 & 16\\
        0 & 0
    \matr,
    \\
    C&=\matl
         4  &  1 & -40 &  -14\\
         -32 & -12 & 4  &  2
    \matr, \quad 
    D=\matl
        0 & 64\\
        64 & 0
    \matr,
\end{align}
and MATLAB's \texttt{tzero} routine yields ${\rm izeros} 
        \left[
        \begin{array}{c|c}
            A & B \\
            \hline
            C & D
        \end{array} 
        \right] = \{-0.77 \pm 1.38i,
  -2.23 \pm 2.16i
  \}.$

The invariant zero decomposition requires $D=0$ and Theorem \ref{theorem:zeros_MIMO} is valid only if $D= 0.$
Therefore, to apply Theorem \ref{theorem:zeros_MIMO}, we consider the system
\begin{align}
    H(s)
        =
            \frac{G(s)}{s+16}.
\end{align}
The relative degrees of each output of $H(s)$ is 1, that is, $\rho_1=\rho_2=1,$ and thus $\rho=2.$
A realization of $H(s)$ is
\begin{align}
    A&=\matl
       -28 &  -14 & -4 & 0 & 0 & 0\\
       16 & 0 & 0 & 0 & 0 & 0\\
       0 & 8 & 0 & 0 & 0 & 0\\
       0 & 0 & 0 & -28 & -14 & -4\\
       0 & 0 & 0 & 16 & 0 & 0\\
       0 & 0 & 0 & 0 & 8 & 0
    \matr, \quad 
    B=\matl
        1 & 0\\
        0 & 0\\
        0 & 0\\
        0 & 1\\
        0 & 0\\
        0 & 0
    \matr,
    \\
    C&=\matl
         0  &  4 & 1 & 64 &  8 & 2\\
         64 & 16 & 4  &  0 & 4 & 2
    \matr, \quad 
    D=\matl
        0 & 0\\
        0 & 0
    \matr.
\end{align}
Note that the pole at $s=-16$ is chosen to obtain integer values in the state-space realization of $H(s).$
Letting 
\begin{align}
    B_z=\matl0 & 0 & 1 & 0 & 0 & 0\\
         0 & 1 & 0 & 0 & 0 & 0\\
         0 & 0 & 0 & 0 & 1 & 0\\
         0 & 0 & 0 & 0 & 0 & 1
         \matr,
\end{align}
it follows from \eqref{eq:T_def} that 
    \begin{align}
        T=\matl
            B_z\\
            C_1\\
            C_2
        \matr=\matl
            0 & 0 & 1 & 0 & 0 & 0\\
            0 & 1 & 0 & 0 & 0 & 0\\
            0 & 0 & 0 & 0 & 1 & 0\\
            0 & 0 & 0 & 0 & 0 & 1\\
            0 & 4 & 1 & 64 & 8 & 2\\
            64 & 16 & 4 & 0 & 4 & 2
        \matr,
    \end{align}
    and thus
    \begin{align}
        \SA&=\left[\begin{array}{cccc|cc}
            0 & 8 & 0 & 0 & 0 & 0\\
            -1 & -4 & -1 & -0.5 & 0 & 0.25\\
            -0.25 & -1 & -2 & -0.5 & 0.25 & 0\\
            0 & 0 & 8 & 0 & 0 & 0\\
            \hline
            22 & 96 & -676 & -206 & -26 & 1\\
            -161 & -484 & 104 & 46 & 1 & -24
        \end{array}\right], 
        \nn \\
        \SB
            &=
                \left[
                \begin{array}{cc}
                0 & 0\\
                0 & 0\\
                \hline
                0 & 0\\
                0 & 0\\
                0 & 64\\
                64 & 0
                \end{array}
                \right]\\
        \SC&=\left[\begin{array}{cccc|cc}
            0 & 0 & 0 & 0 & 1 & 0\\
            0 & 0 & 0 & 0 & 0 & 1
        \end{array}\right].
    \end{align}
    Note that 
    \begin{align}
        \SA_{\eta}
            =
                \matl
                    0 & 8 & 0 & 0\\
                    -1 & -4 & -1 & -0.5\\
                    -0.25 & -1 & -2 & -0.5\\
                    0 & 0 & 8 & 0\\
                \matr
    \end{align}
    and thus $\text{mspec}(\SA_{\eta})=\{-2.23 \pm 2.16i,
  -0.77 \pm 1.38i\}$, which confirms Theorem \ref{theorem:zeros_MIMO}.
\end{exmp}

\begin{exmp} 
\textbf{[Nonsquare MIMO system.]}
Consider the system
\begin{align}
    G(s)=\dfrac{1}{s^3 +  2s^2 -s -2}\matl
        s^2 + s -2 &  0 & s^2 -2s + 1\\
        s^2 - 3 s + 2 &  s^2 - 1 & s^2 - 1
    \matr,
\end{align}
Note that $\rho_1=1, \rho_2=1$, and thus $\rho=2$. A realization of $G$, computed with MATLAB's \texttt{ss} routine, is
\begin{align}
    A&=\matl
       0 &  0 & 2 & 0 & 0 & 0\\
       1 & 0 & 1 & 0 & 0 & 0\\
       0 & 1 & -2 & 0 & 0 & 0\\
       0 & 0 & 0 & 0 & 0 & 2\\
       0 & 0 & 0 & 1 & 0 & 1\\
       0 & 0 & 0 & 0 & 1 & -2
    \matr, \quad 
    B=\matl
        -2 & 0 & 1\\
        1 & 0 & -2\\
        1 & 0 & 1\\
        2 & -1 & -1\\
        -3 & 0 & 0\\
        1 & 1 & 1
    \matr,
    \\
    C&=\matl
         0  &  0 & 1 &  0 & 0 & 0\\
         0  &  0 & 0 &  0 & 0 & 1
    \matr, \quad 
    D=\matl
        0 & 0 & 0\\
        0 & 0 & 0
    \matr,
\end{align}
and MATLAB's \texttt{tzero} routine yields ${\rm izeros} 
        \left[
        \begin{array}{c|c}
            A & B \\
            \hline
            C & D
        \end{array} 
        \right] = \{1,1
  \}.$

Note that Theorem \ref{theorem:zeros_MIMO} is proved only for square systems.
To apply Theorem \ref{theorem:zeros_MIMO} to a nonsquare MIMO system, we square the system 
and use the invariant zero decomposition to compute the invariant zeros of the squared system. 
The invariant zeros of the original nonsquare system can be computed by either checking the rank of the Rosenbrock system matrix or by squaring the nonsquare system a few times and collecting the common invariant zeros.

In this example, we square the system by augmenting the output matrix $C$ with additional rows to match the number of columns of $B$. 
Note that $D$ is extended with zero rows. 
In particular, $C$ is augmented as follows. 
\begin{align}
    C_{\rm sq 1}
        &=
            \matl
               0  &  0 & 1 &  0 & 0 & 0\\
               0  &  0 & 0 &  0 & 0 & 1\\
               2  &  1 &  3 & 4 & 7 & 3
            \matr, 
    \nn\\
    C_{\rm sq 2}
        &=
            \matl
               0  &  0 & 1 &  0 & 0 & 0\\
               0  &  0 & 0 &  0 & 0 & 1\\
               1  &  4 & 2 &  4 & 0 & 3
            \matr, 
    \nn \\
    D_{\rm sq }
        &=
            0_{3\times 3}.
\end{align}
Note that for both squared systems $(A,B,C_{\rm sq1},D_{\rm sq} )$ and $(A,B,C_{\rm sq2},D_{\rm sq} )$, 
$\rho_1= \rho_2= \rho_3=1$, and thus $\rho=3$ in both cases.
Letting 
{\footnotesize
\begin{align}
    B_z&=\matl
        0.6427  &  0.4239  &  0.2051    &0.4063  &  0.1875  &  0.4063\\
        -0.5321  & -0.0333  &  0.4654    &0.1671  &  0.6658  &  0.1671\\
        -0.1106  & -0.3906  & -0.6706    &0.4267  &  0.1467  &  0.4267
         \matr,
\end{align}
}
it follows from \eqref{eq:T_def} that 
    {\footnotesize
    \begin{align}
        T_{\rm sq 1}
            =&
                \matl
                    B_z\\
                    C_{\rm sq 1}
                \matr
            \nn\\
            =&
                \matl
                    0.6427  &  0.4239  &  0.2051 &0.4063  &  0.1875  &  0.4063\\
                    -0.5321  & -0.0333  &  0.4654  &  0.1671  &  0.6658 &0.1671\\
                    -0.1106  & -0.3906 &  -0.6706  &  0.4267  &  0.1467  &  0.4267\\
                    0  & 0  &  1  & 0  & 0 & 0\\
                    0  & 0  &  0  & 0  & 0 & 1\\
                    2  & 1  &  3  & 4  & 7 & 3
                \matr,
        \\
        T_{\rm sq 2}
            =&
                \matl
                    B_z\\
                    C_{\rm sq 2}
                \matr
            \nn\\
            =&
                \matl
                    0.6427  & 0.4239  & 0.2051   &0.4063  & 0.1875  & 0.4063\\
                    -0.5321  & -0.0333  & 0.4654 &0.1671  & 0.6658  & 0.1671\\
                    -0.1106  & -0.3906 & -0.6706 &0.4267  & 0.1467  & 0.4267\\
                    0  & 0  & 1 & 0 & 0 & 0\\
                    0  & 0  & 0 & 0 & 0 & 1\\
                    1  & 4  & 2 & 4 & 0 & 3
                \matr.
    \end{align}
    }
Note that 
    \begin{align}
        \SA_{\eta1}
            &=
                \matl
                   -0.1009 &  -0.4415   &-0.3170\\
                    2.5100  &  2.0067 &0.7227\\
                  -1.4090 &  -0.5651    &0.5943
                \matr,
            \nn \\
        \SA_{\eta2}
            &=
                \matl
                   1.0241  &  0.3400   &-0.0829\\
                  -0.0550  &  0.2248    &0.1889\\
                  0.0309  &  0.4352    &0.8939
                \matr
    \end{align}
    and thus $\text{mspec}(\SA_{\eta1})=\{0.5,1,1\}$, and $\text{mspec}(\SA_{\eta2})=\{0.14,1,1\}$. The common eigenvalues for $\SA_{\eta1}$ and $\SA_{\eta2}$ are $\{1,1\}$, which are the same as the values obtained by MATLAB's \texttt{tzero} routine.
\end{exmp}

\section{Conclusion}
\label{sec:conclusion}

In this paper, we introduced the invariant zero form of a MIMO linear system, which is a realization that isolates the invariant zeros of a MIMO linear system into a partition of the transformed dynamics matrix.
%
We constructed a change of basis matrix that transforms a given realization of a MIMO linear system into its invariant zero form. 
Furthermore, we showed that, in the invariant zero form, the transformed dynamics, input, and output matrices have sparse structure.
Finally, we showed that, in the case of a square MIMO system, the invariant zeros of the system are isolated in a partition of the dynamics matrix in the invariant zero form and can be computed by solving an eigenvalue problem.  
Numerical examples are presented to confirm the paper's main result and extend its application to the nonsquare cases. 
%


%
%
%

\bibliographystyle{elsarticle-num}
\bibliography{b}

\begin{thebibliography}{10}
\expandafter\ifx\csname url\endcsname\relax
  \def\url#1{\texttt{#1}}\fi
\expandafter\ifx\csname urlprefix\endcsname\relax\def\urlprefix{URL }\fi
\expandafter\ifx\csname href\endcsname\relax
  \def\href#1#2{#2} \def\path#1{#1}\fi

\bibitem{macfarlane1976poles}
A.~MacFarlane, N.~Karcanias, Poles and zeros of linear multivariable systems: a survey of the algebraic, geometric and complex-variable theory, International Journal of Control 24~(1) (1976) 33--74.

\bibitem{desoer1974zeros}
C.~Desoer, J.~Schulman, Zeros and poles of matrix transfer functions and their dynamical interpretation, IEEE Transactions on Circuits and Systems 21~(1) (1974) 3--8.

\bibitem{tokarzewski2006finite}
J.~Tokarzewski, Finite Zeros in Discrete Time Control Systems, Vol. 338, Springer, Berlin, 2006.

\bibitem{hoagg2007nonminimum}
J.~B. Hoagg, D.~S. Bernstein, Nonminimum-phase zeros-much to do about nothing-classical control-revisited part ii, IEEE Control Systems Magazine 27~(3) (2007) 45--57.

\bibitem{havre2001achievable}
K.~Havre, S.~Skogestad, Achievable performance of multivariable systems with unstable zeros and poles, International Journal of Control 74~(11) (2001) 1131--1139.

\bibitem{wonham1970decoupling}
W.~M. Wonham, A.~S. Morse, Decoupling and pole assignment in linear multivariable systems: a geometric approach, SIAM Journal on Control 8~(1) (1970) 1--18.

\bibitem{owens1977invariant}
D.~Owens, Invariant zeros of multivariable systems: A geometric analysis, International Journal of Control 26~(4) (1977) 537--548.

\bibitem{davison1974properties}
E.~Davison, S.~Wang, Properties and calculation of transmission zeros of linear multivariable systems, Automatica 10~(6) (1974) 643--658.

\bibitem{horan1980decoupling}
R.~Horan, On the decoupling zeros and poles of a system, IEEE Transactions on Automatic Control 25~(3) (1980) 517--521.

\bibitem{patel1986blocking}
R.~Patel, On blocking zeros in linear multivariable systems, IEEE transactions on automatic control 31~(3) (1986) 239--241.

\bibitem{laub1978calculation}
A.~J. Laub, B.~Moore, Calculation of transmission zeros using {QZ} techniques, Automatica 14~(6) (1978) 557--566.

\bibitem{tokarzewski2011invariant}
J.~Tokarzewski, Invariant zeros of linear singular systems via the generalized eigenvalue problem, in: Future Intelligent Information Systems, Springer, 2011, pp. 263--270.

\bibitem{tokarzewski2009zeros}
J.~Tokarzewski, Zeros in linear control systems and generalized eigenvalue problem, PRZEGLAD ELEKTROTECHNICZNY 85~(10) (2009) 129--132.

\bibitem{bernstein2009matrix}
D.~S. Bernstein, Matrix mathematics, in: Matrix Mathematics, Princeton university press, Ann Arbor, 2009.

\bibitem{emami1982computation}
A.~Emami-Naeini, P.~Van~Dooren, Computation of zeros of linear multivariable systems, Automatica 18~(4) (1982) 415--430.

\bibitem{misra1994computation}
P.~Misra, P.~Van~Dooren, A.~Varga, Computation of structural invariants of generalized state-space systems, Automatica 30~(12) (1994) 1921--1936.

\bibitem{misra1989computation}
P.~Misra, R.~V. Patel, Computation of minimal-order realizations of generalized state-space systems, Circuits, Systems and Signal Processing 8~(1) (1989) 49--70.

\bibitem{misra1989minimal}
P.~Misra, R.~V. Patel, Minimal order generalized state space representation of singular systems, in: 1989 American Control Conference, IEEE, 1989, pp. 2140--2145.

\bibitem{davison1978algorithm}
E.~Davison, S.~Wang, An algorithm for the calculation of transmission zeros of the system {(C, A, B, D)} using high gain output feedback, IEEE Transactions on Automatic Control 23~(4) (1978) 738--741.

\bibitem{garbow1977matrix}
B.~S. Garbow, J.~M. Boyle, J.~J. Dongarra, C.~B. Moler, Matrix eigensystem routines—EISPACK guide extension, Springer, Berlin, 1977.

\bibitem{aling1984nine}
H.~Aling, J.~M. Schumacher, Nine-fold canonical decomposition for linear systems., International Journal of control 39~(4) (1984) 779--806.

\bibitem{basile1994controlled}
G.~Basile, G.~Marro, J.~Schumacher, Controlled and conditioned invariants in linear system theory, IEEE Transactions on Automatic Control 39~(1) (1994) 250--250.

\bibitem{morris2010invariant}
K.~Morris, R.~Rebarber, Invariant zeros of siso infinite-dimensional systems, International journal of control 83~(12) (2010) 2573--2579.

\bibitem{isidori1985nonlinear}
A.~Isidori, Nonlinear control systems: an introduction, Springer, Italy, 1985.

\bibitem{zhao2024weak}
L.~Zhao, Z.-Y. Li, B.~Zhou, W.~Michiels, On the weak infinite zero structures and relative degrees for linear time-delay systems, SIAM Journal on Matrix Analysis and Applications 45~(4) (2024) 1873--1901.

\bibitem{zhou2024invertibility}
B.~Zhou, On the invertibility indices and the morse normal form for linear multivariable systems, International Journal of Control 97~(6) (2024) 1198--1209.

\bibitem{portella2024circumventing}
J.~M. Portella~Delgado, A.~Goel, Circumventing unstable zero dynamics in input-output linearization of longitudinal flight dynamics, in: AIAA SCITECH 2024 Forum, 2024, p. 0321.

\bibitem{zhou2023relative}
B.~Zhou, On the relative degree and normal forms of linear systems by output transformation with applications to tracking, Automatica 148 (2023) 110800.

\bibitem{korovin2007canonical}
S.~Korovin, A.~Il’in, V.~Fomichev, A canonical form of vector control systems, in: Doklady Mathematics, Vol.~75, Springer, 2007, pp. 467--471.

\bibitem{portella2024computing}
J.~M.~P. Delgado, A.~Goel, Computing invariant zeros of a linear system using state-space realization, in: 2024 American Control Conference (ACC), IEEE, 2024, pp. 2746--2751.

\bibitem{Khalil:1173048}
H.~K. Khalil, Nonlinear systems; 3rd ed., Prentice-Hall, Upper Saddle River, NJ, 2002.

\bibitem{isidori1988nonlinear}
A.~Isidori, C.~H. Moog, On the nonlinear equivalent of the notion of transmission zeros, in: Modelling and Adaptive Control, Springer, 1988, pp. 146--158.

\bibitem{tian2004rank}
Y.~Tian, Rank equalities for block matrices and their moore-penrose inverses, Houston J. Math 30~(4) (2004) 483--510.

\bibitem{buaphim2018some}
N.~Buaphim, K.~Onsaard, P.~So-ngoen, T.~Rungratgasame, Some reviews on ranks of upper triangular block matrices over a skew field, in: International Mathematical Forum, Vol.~13, 2018, pp. 323--335.

\end{thebibliography}

\appendix

\renewcommand{\thesection}{A}
\section{Algorithms to Compute Zeros of a System}
\renewcommand{\thesection}{A}
\label{appndx:Algorithms}

This appendix presents the pseudo-code to compute the invariant zeros of a MIMO linear system using the GAZERO algorithm presented in \cite{basile1994controlled} and the invariant zero decomposition algorithm presented in this paper.

\begin{algorithm}
\DontPrintSemicolon
  
  \KwInput{$A,B,C,D$ \tcp*{System matrices }}
  \KwOutput{$z$ \tcp*{Column vector containing the invariant zeros}}
    
  $\text{kerC} \leftarrow \text{ker}(\text{C})$ \tcp*{orthogonal basis for the nullspace of C}
  
  $\text{imB} \leftarrow \text{ima}(\text{B})$\tcp*{orthogonal basis for the image of B}
  
  $\mathcal{V}^* \leftarrow \text{mainco}(\text{A},\text{imB},\text{kerC})$ \tcp*{ maximal locus of controlled state trajectories such that the output is identically zero}

      \If{$\mathcal{V}^* = \emptyset$}
      {
        $X_{22} \leftarrow  0$ \tcp*{ denotes the invariant zero structure}
        $z \leftarrow \emptyset$, \Return \tcp*{ invariant zeros of the system}
      }
  $\mathcal{S}^* \leftarrow \text{miinco}(\text{A},\text{kerC},\text{imB})$ \tcp*{  the set of all states reachable from the origin through trajectories having all states but the ones belonging to the nullspace of C}

  $\mathcal{R}^* \leftarrow \text{ints}(\mathcal{V}^*,\mathcal{S}^*)$ \tcp*{gives an orthonormal basis for the subspace image $\mathcal{V}^*$ intersected by image $\mathcal{S}^*$}

    \If{$\mathcal{R}^* = \emptyset$}
      {
        $V_1 \leftarrow \mathcal{V}^*$
      }
    \Else
    {
        $V_1 \leftarrow \text{ima}([\mathcal{R}^*,\mathcal{V}^*],0)$ \tcp*{ Suitable orthonormal basis for matrix $V_1$, needed for transformation, the $0$ value is to not allow permutations}
    }

  
  $X_{22} \leftarrow \text{pinv}([ V_1 \,\text{B}])\text{A}V_1$ \tcp*{ denotes the $(nV-nR) \times (nV-nR)$ submatrix extracted from the bottom-right corner of matrix $X$, whose image is denoted by the $(A,B)$-controlled invariant $\mathcal{V}$ subspace.}

  $z \leftarrow \text{eig}(X_{22})$ \tcp*{ Invariant zeros of the system.}

  \Return
  
\caption{GAZERO Algorithm}
\end{algorithm}

\begin{algorithm}
\DontPrintSemicolon
  
  \KwInput{$A,B,C,D$ \tcp*{System matrices }}
  \KwOutput{$z$ \tcp*{Invariant zeros}}
  $\overline{B} \leftarrow \text{null}(B^{\rm T})$\\
  $\rho \leftarrow \text{RelativeDegree}(System)$ 
  
  \While{$j \leq l_y$}
  {
  \While{$i \leq l_y$}
   {
   		$\overline{C}_j \leftarrow \text{stack}(C_iA^{i-1})[\text{bottom}]$\\
   }
   $\overline{C} \leftarrow \text{stack}(C_j)[\text{bottom}]$\\
  }
  \If { $\exists (l_x - \rho) \text{\rm independent rows of } \overline{B} \text{\rm with respect to } \overline{C}$}
  {
    $B_z \leftarrow (l_x - \rho)\text{  independent rows of } \overline{B}$\\
    $T \leftarrow [B_z, \overline{C}]^{\rm T}$
  }
  \Else
  {
    $T$ \text{ is not invertible}, \Return
  }
  $\mathcal{A} \leftarrow T A T^{-1}$\\
  $\mathcal{A}_{\eta} \leftarrow (l_z \times l_z) \text{ upper-left block of }\mathcal{A}$\\
  $z \leftarrow \text{eig}(\mathcal{A}_{\eta})$ \tcp*{Invariant zeros}
  \Return
\caption{Invariant Zero Decomposition Algorithm}
\end{algorithm}

\renewcommand{\thesection}{B}

\section{Dynamic Extension of Linear System}
\renewcommand{\theequation}{\thesection.\arabic{equation}}
\setcounter{equation}{0}
\renewcommand{\thesection}{B}
\label{appndx:feedforward}
Consider the system 
\begin{align}
    \dot x &= A x + B u ,
    \label{eq:appndx_xdot}
    \\
    y &= C x + D u,
    \label{eq:appndx_y}
\end{align}
where 
$x \in \BBR^{l_x}.$
Define
\begin{align}
    v \isdef \dot u + \alpha I_{l_u} u.
\end{align}
Defining $\overline x = \matl x \\ u \matr, $ it follows that 
\begin{align}
    \dot {\overline x }
        &=
            \overline A
            \overline x
            +
            \overline B
            v
    \label{eq:appndx_xbardot}
            \\
    y
        &=
            \overline C
            \overline x,
    \label{eq:appndx_ybar}
\end{align}
where 
\begin{align}
    \overline A 
        \isdef 
            \matl 
                A & B \\
                0 & -\alpha I_{l_u}
            \matr 
    , \quad
    \overline B 
        \isdef 
            \matl 
                0_{l_x \times l_u} \\
                I_{l_u}
            \matr 
    , \quad 
    \overline C
        \isdef 
            \matl 
                C  & D
            \matr .
\end{align}
Note that \eqref{eq:appndx_xbardot}-\eqref{eq:appndx_ybar} is a \textit{dynamic extension} of \eqref{eq:appndx_xdot}-\eqref{eq:appndx_y}.

\begin{theorem}
    \label{theo:GH_zeros}
    Consider the system \eqref{eq:appndx_xdot}-\eqref{eq:appndx_y} and its extension \eqref{eq:appndx_xbardot}-\eqref{eq:appndx_ybar}. 
    Then, the invariant zeros of $(A,B,C,D)$ are the same as the invariant zeros of $(\overline{A}, \overline{B}, \overline{C}, 0_{l_y \times l_u} ).$
\end{theorem}
\begin{proof}
    Note that 
\begin{align}
    \overline \chi (s) 
        =
            \matl
                sI-\overline A & \overline B\\
                \overline C & -\overline D
            \matr
        =
            \matl
                sI- A & -B & 0_{l_x \times l_u}\\
                0 & (s-\alpha) I_{l_u} & I_{l_u} \\ 
                C &  D & 0_{l_y \times l_u}
            \matr. \nn
\end{align}
Since the rank of a matrix remains unchanged by a row exchange operation, it follows that 
\begin{align}
    \rank \  \overline \chi (s) 
        =
            \rank \  
            \matl
                sI- A & B & 0_{l_x \times l_u}\\ 
                C &  -D & 0_{l_y \times l_u} \\
                0 & (s-\alpha) I_{l_u} & I_{l_u}
            \matr. \nn
\end{align}
Next, it follows from Lemma \ref{lem:rank_tian} that 
\begin{align}
    \rank \  \matl \SC & 0 \\ \SA & \SB \matr
        =&
            \rank \  \SB + \rank \  \SC
            \nn\\
            &+ \rank \  \left[(I - \SB \SB^+) \SA (I - \SC^+ \SC) \right].
            \nn
\end{align}
Setting $\SC = \chi(s)$ and $\SB = I_{l_u}$ yields
\begin{align}
    \rank \  \overline \chi (s) 
        =
            \rank \  \chi (s)
            +
            \rank \  I_{l_u}
        =
            \rank \  \chi (s)
            +
            l_u, \nn
\end{align}
thus implying that the $\rank \  \overline{\chi}(s)$ drops if and only if $\rank \  \chi(s)$ drops, which completes the proof.   
\end{proof}


\renewcommand{\thesection}{C}
\setcounter{equation}{0}
\section{Linear Algebra Facts}
\renewcommand{\thesection}{C}
\label{appndx:LinearAlgebraFacts}
This appendix reviews several linear algebra facts used to prove various results in the paper. 


The following Lemma is given in \cite{tian2004rank}, \cite{buaphim2018some}, and is used in the proof of Proposition \ref{prop:delta_n_rank} and Theorem \ref{theo:GH_zeros}
\begin{lemma} \label{lem:rank_tian}
     Let 
    $A \in \BBC^{m \times n},$
    $B \in \BBC^{m \times k},$ and
    $C \in \BBC^{l \times n}.$
    Then, 
    \begin{align}
        \rank \matl A & B \\ C & 0 \matr
            &=
                \rank B + \rank C
                \nn\\ 
                &+ \rank \left[(I - B B^+) A (I - C^+ C) \right],
    \end{align}    
and
    \begin{align} 
        \rank \matl
            A\\
            C
        \matr 
            &=
                \rank (A) + \rank (C - CA^{{+}}A)
            \\
            &=
                \rank (C) + \rank (A - AC^{{+}}C)
            \\
            &=
                \rank(A) + \rank(C) - \dim\left(\mathcal{R}(A^{\rm T}) \cap \mathcal{R}(C^{\rm T})\right)
            \label{eq:rankAC}
            \\
            &\geq \max (\rank(A), \rank(C)).
    \end{align}
\end{lemma}


The following results are used in the proof of Theorem \ref{theorem:zeros_MIMO}.
\begin{proposition}
    \label{prop:delta_n_rank}
    Let $n \in \BBN.$
    For all $s \in \BBC,$ the rank of 
\begin{align}
    \Delta_{n}(s)
        &\isdef
            \matl
                -1 & 0 & \cdots & 0 & 0\\
                s & -1 & \cdots & 0 & 0\\
                \vdots & \vdots & \ddots & \vdots & \vdots\\
                0 & 0 & \cdots & -1 & 0\\
                0 & 0 & \cdots & s & -1
            \matr
    \in \mathbb{R}^{n\times n} 
    \label{def:delta_n}
\end{align}
is $n.$\\
\end{proposition}

\begin{proof}
Note that for all $s\in \BBC,$ 
$\det\left(\Delta_n(s)\right)=\pm1,$ which implies that $\rank (\Delta_n(s))=n.$    
\end{proof}

\begin{proposition}
    \label{prop:A_Aplus_zeros}
    Let $A \in \mathbb{R}^{n\times m}.$
    Let $A^+$ be the Moore-Penrose psuedo-inverse of $A.$
    If the $i$th row of $A$ is zero, then the $i$th row and $i$th column of $A A^{+}$ is zero. 
    %
    %
    %
\end{proposition}

\begin{proof}
    Let the $i$th row of $A$ be zero. 
    Then, there exists a permutation matrix $P \in \BBR^{n \times n}$ such that
    \begin{align}
        B
            \isdef 
                PA
            =
            \matl
                0_{1 \times m}\\
                \overline A
            \matr, \nn 
    \end{align}
    where $\overline A$ contains all but the $i$th row of $A.$
    Note that $A = P^{\rmT} B.$
    Then, 
    \begin{align}
        BB^{{+}} 
            =
                \matl 
                0_{1\times m}\\
                \overline{A}
                \matr\matl
                0_{m \times 1} & \overline{A}^{{+}}
                \matr=\matl
                0 & 0_{1 \times (n-1)}\\
                0_{(n-1) \times 1} &  \overline{A}\,\, \overline{A}^+
                \matr.
        \nn 
    \end{align}
    %
    Finally, note that $A A^+ = P^{\rmT}B B^+ P $ where the pre-multiplication by $P^\rmT$ and post-multiplication by $P$ in $P^\rmT B B^+ P $ permutes the first zero row and the first zero column to the appropriate position.
    \end{proof}

\begin{proposition}
    \label{prop:Aplus_A_zeros_column}
    Let $A^{n\times m}$.
    If the $i$th column of $A$ is zero, then the $i$th row and $i$th column of $A^{+} A$ is zero.
\end{proposition}

\begin{proof}
    The proof is similar to the proof of proposition \ref{prop:A_Aplus_zeros}.
\end{proof}

\clearpage

\end{document}